\documentclass{notices}
\usepackage{amsfonts,amssymb,amsmath,amscd,graphicx}

\def\complex{{\mathbb C}}

    \def\amano{1}
   \def\banjai{2}
\def\bornemann{3}
   \def\crowdy{4}
     \def\dekp{5}
      \def\dlg{6}
 \def\driscoll{7}
  \def\chebfun{8}
       \def\dt{9}
     \def\crdt{10}
    \def\gaier{11}
    \def\gopal{12}
      \def\rep{13}
  \def\henrici{14}
   \def\zipper{15}
   \def\nasser{16}
 \def\odonnell{17}
     \def\papa{18}
 \def\schiffer{19}
  \def\packing{20}
   \def\sisc80{21}
   \def\series{22}
      \def\ncm{23}
  \def\wegmann{24}

\title{Numerical Conformal Mapping}

\author{
Lloyd N. Trefethen
\affil{Professor of Applied Mathematics in Residence,
School of Engineering and Applied Sciences, Harvard University}
}

\begin{document}

\maketitle

Conformal mapping may be the best-known topic in complex analysis.
Any simply connected nonempty domain $\Omega$ in the complex plane
$\complex$ (assuming $\Omega\ne \complex$) can be mapped bijectively
to the unit disk by an analytic function with nonvanishing
derivative, as in Figure 1.  If $\Omega$ is doubly-connected, it can
be mapped to a circular annulus $1<|z|<R$ for some $R$, called the
{\em conformal modulus,} which is uniquely determined by $\Omega$,
as in Figure 2.  If $\Omega$ has connectivity higher than $2$, it
can be mapped onto various canonical domains such as a disk with
exclusions in the form of slits or smaller disks, as in Figure 3.

\begin{figure}[h]
\noindent\kern 8pt\includegraphics [trim=30 150 100 45, clip, width=3.0in]{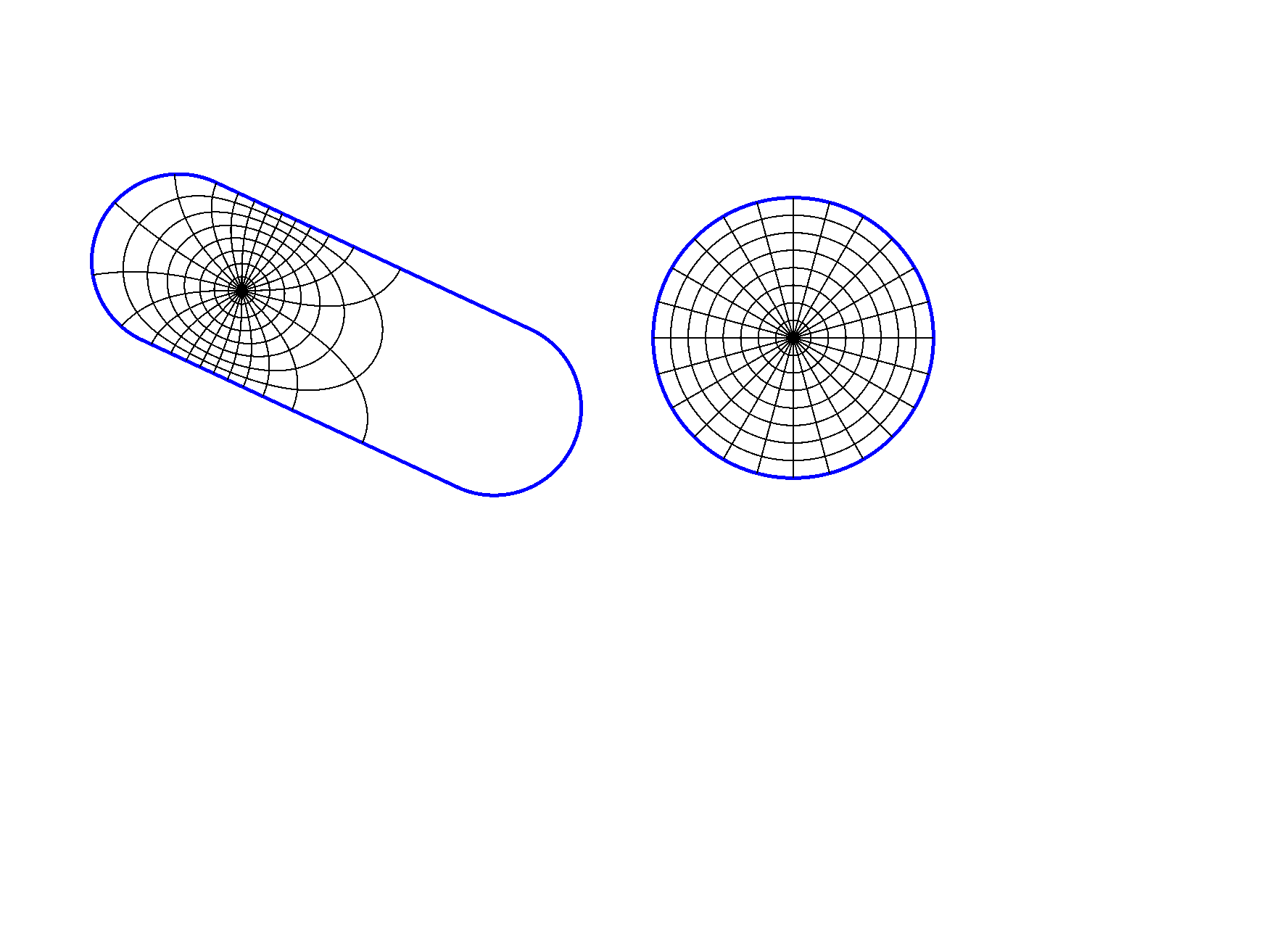}
\vspace*{-3pt}
\caption{A simply-connected conformal map onto the unit disk.
The most common target domains
are a disk, a half-plane, or a rectangle.}
\end{figure}

\begin{figure}[h]
\noindent\kern 5pt\includegraphics [trim=30 145 50 45, clip, width=3.3in]{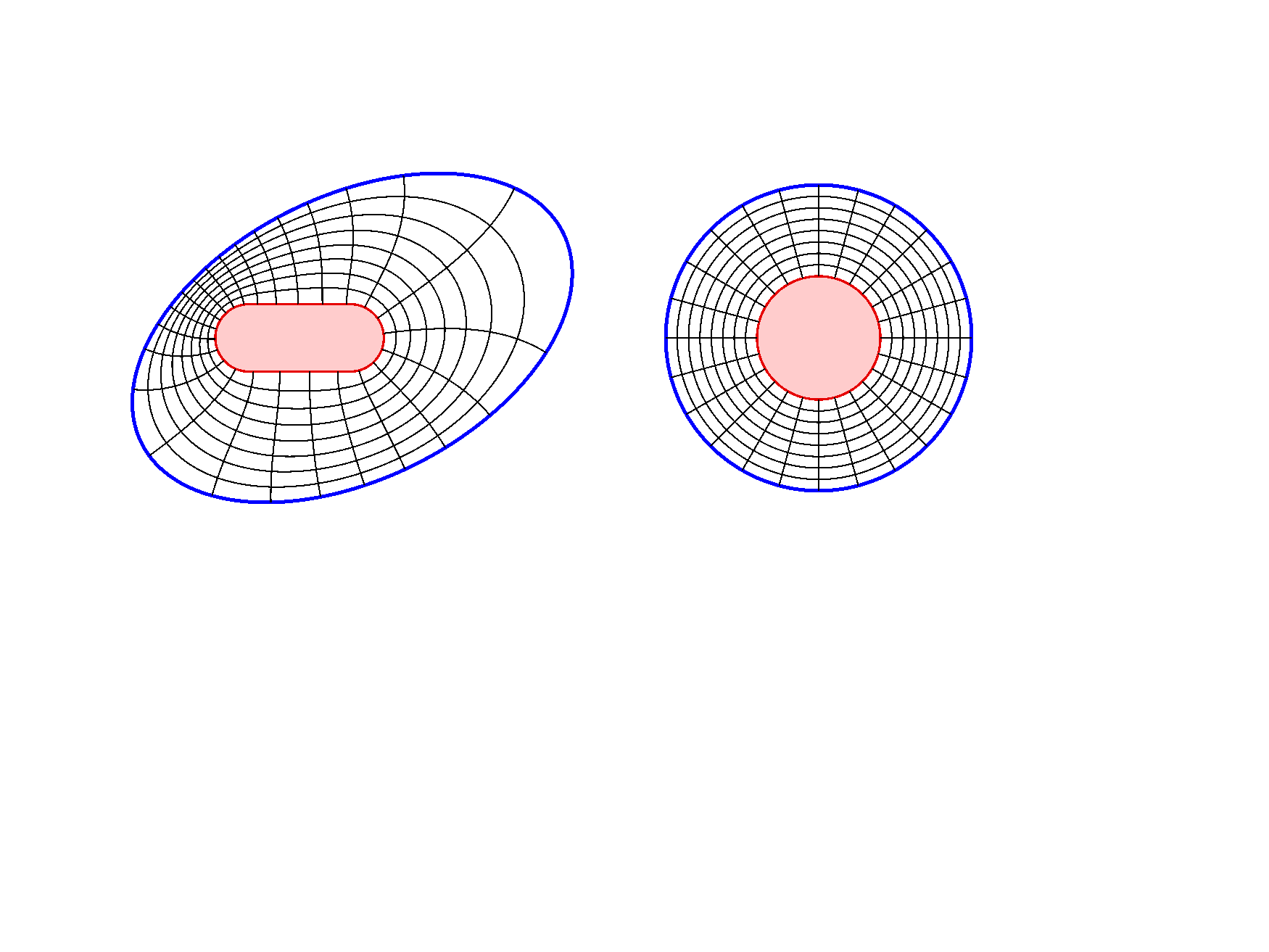}
\vspace*{-15pt}
\caption{A doubly-connected conformal map onto a circular annulus,
which is the most common target domain.}
\end{figure}

\begin{figure}[h]
\noindent\kern 6pt\includegraphics [trim=30 145 50 45, clip, width=3.7in]{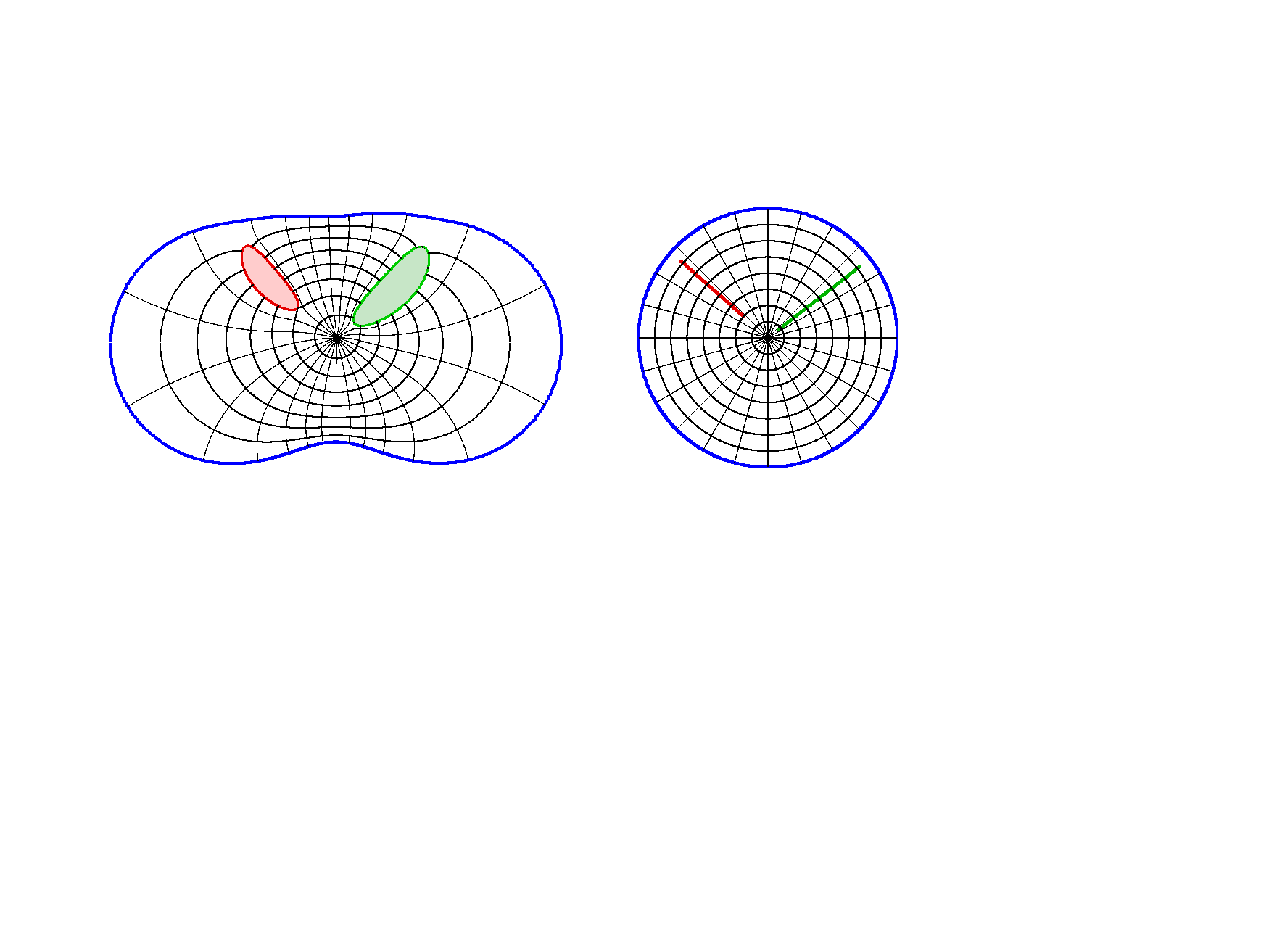}
\vspace*{-23pt}
\caption{A triply-connected conformal map onto a disk with radial
slits.  Other standard target domains involve exclusions in the form of disks or
circular arcs.}
\end{figure}

Most conformal maps cannot be found analytically, so when
computers began to appear in the 1950\kern .5pt s and 1960\kern
.5pt s, the field of {\em numerical conformal mapping} was born.
Many methods involve discretizations of integral equations, such
as those of Gerschgorin, Lichtenstein, Kerzmann-Stein-Trummer, Symm,
and Theodorsen.  Dieter Gaier of the University of Giessen published
an important monograph in those early years, which unfortunately was
never translated to English [\gaier].  A later reference is the book
of Henrici [\henrici], and there is a major survey paper by Wegmann
[\wegmann].

Conformal mapping of polygons has always been a conspicuous special
case, thanks to the Schwarz-Christoffel transformation, and this
became readily available for applications with the appearance
of Driscoll's SC Toolbox in Matlab in the 1990\kern .5pt s
[\driscoll,\dt,\sisc80].  To this day, this Toolbox remains the most
widely-used software for conformal mapping, and Figure 4 illustrates
that its power extends to nontrivial geometries [\crdt].  This figure
shows the mapping of a {\em quadrilateral,} namely a domain with four
distinguished boundary points, onto a rectangle, whose aspect ratio
$\mu \approx 18.20539$
(also called the conformal modulus) is uniquely determined [\papa].

\begin{figure}[h]
\vspace*{7pt}
\noindent\kern 27pt\includegraphics [trim=0 0 0 0, clip, width=2.4in]{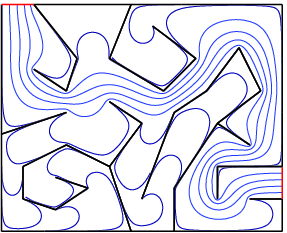}
\caption{Schwarz-Christoffel map of ``Emma's maze'' computed with the SC Toolbox; image
taken from ~[\driscoll].  The problem domain, defined by a rectangle with
crooked slits, is known as a quadrilateral because of the four
distinguished vertices (where red meets black).  The
canonical domain, not shown, is a rectangle.}
\end{figure}

Many other conformal mapping methods have been developed too.
For multiply-connected domains, for example, there are methods
based on generalized Schwarz-Christoffel mapping [\dekp], the
Schottky-Klein prime function [\crowdy,\dlg], and an integral equation
related to the Neumann kernel
[\nasser].  In a field as old as this, there
are numerous further methods that have been explored, including
``circle packing'' (introduced by Thurston) [\packing], the
``zipper algorithm'' [\zipper], the ``charge simulation method''
(a version of the method of fundamental solutions) [\amano],
rational approximation [\gopal,\ncm], and many varieties of series and
iterations, sometimes accelerated by the Fast Multipole Method [\banjai,\odonnell].
Commands {\tt conformal} and {\tt conformal2} for smooth simply
and doubly connected conformal mapping can be found in
Chebfun [\chebfun].  Figure~5 shows another map of a quadrilateral,
now one with curved sides, computed by the ``lightning''
rational function method [\gopal].

\begin{figure}[h]
\vspace*{7pt}
\noindent\kern 6pt\includegraphics [trim=60 110 50 100, clip, width=3in]{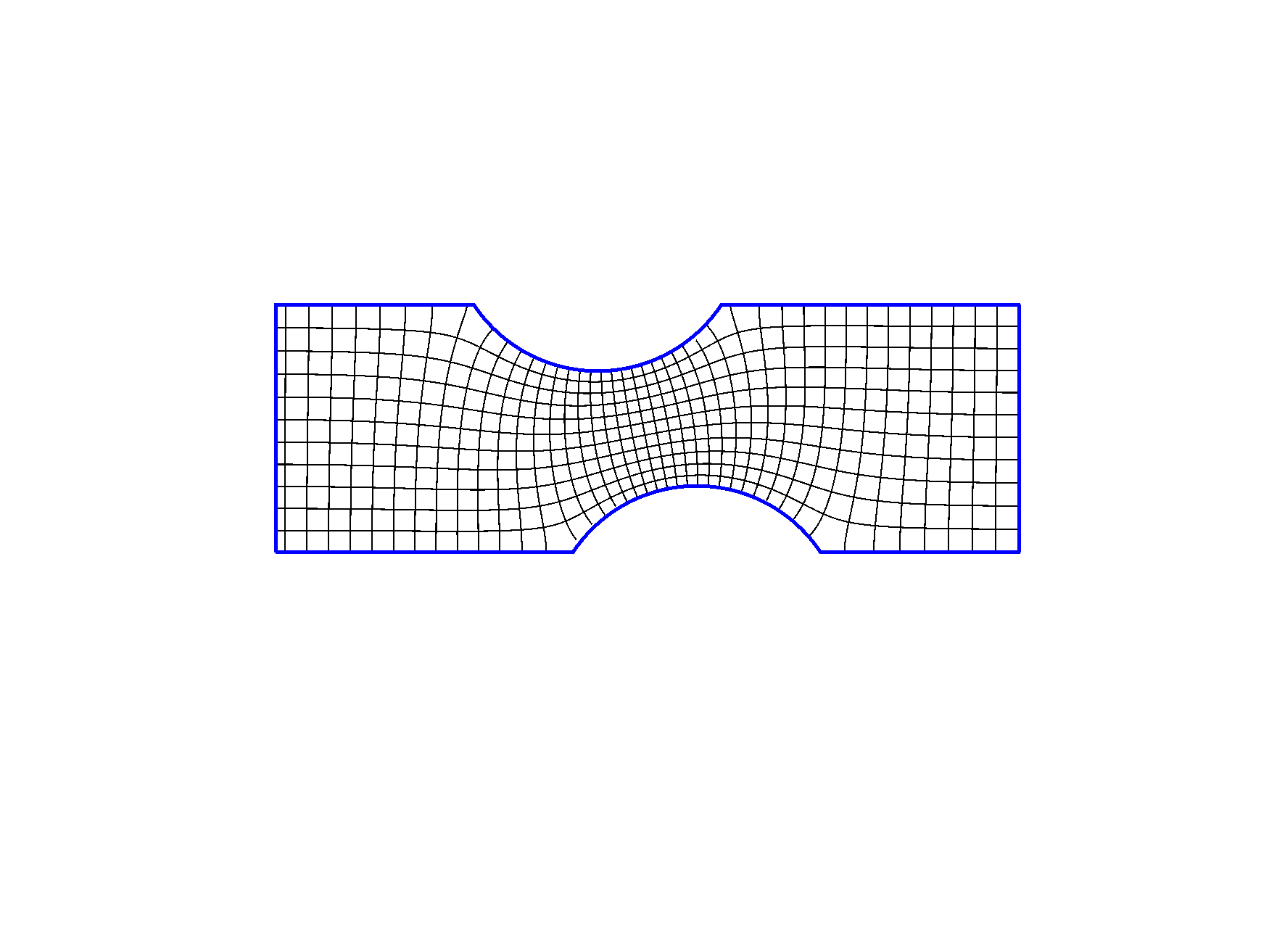}
\vspace*{-18pt}
\caption{Map of another quadrilateral onto a rectangle computed with rational
functions [\gopal,\ncm].}
\end{figure}

Traditionally, a fundamental distinction in numerical conformal
mapping was between methods mapping from problem domain to
canonical domain and those going in the inverse direction, from
canonical domain to problem domain.  In the former case, the map
can be reduced to a Laplace problem, so the integral equations
of Lichtenstein and Symm are linear, for example, whereas that of
Theodorsen in the other direction is nonlinear.  However, the recent
appearance of fast numerical methods for rational approximation has
diminished the distinction between the forward and
inverse maps.  Once one has the {\em boundary correspondence function,}
the homeomorphism between the domain and range boundaries (provided this
makes sense, as is always the case
for Jordan domains), it is an easy matter to use it to compute efficient and
accurate rational approximations in both directions [\rep].
For example, the conformal map of an irregular hexagon
to the unit disk, and its inverse map, can each typically be approximated
to 6-digit accuracy by rational functions of degrees of the order of 50 to
100.\ \ Evaluation of such functions takes just microseconds per point.

As just mentioned, conformal mapping is a special case of a Laplace problem.
For a simply-connected domain $\Omega$ bounded
by a Jordan curve $\Gamma$ enclosing the origin,
following Theorem 16.5a of [\henrici] and p.~253 of [\schiffer], 
suppose we seek the unique conformal map $f$ onto the unit disk with
$f(0) = 0$ and $f'(0)>0$.  
Then $g(z) = \log(f(z)/z)$ is a nonzero
analytic function on $\Omega$ that is continuous on $\overline{\Omega}$ and
has real part $-\log |z|$ for $z\in \Gamma$.  If we write
$g(z) = u(z) + iv(z)$, where $u$ and $v$ are real harmonic functions,
then $u$ is the solution of the Dirichlet problem
\begin{equation}
\Delta u = 0\, ; \quad u(z) = -\log|z|, ~ z\in \Gamma,
\end{equation}
and $v$ is its harmonic conjugate in $\Omega$ with $v(0) = 0$.
Combining these elements, we see that $f$ is given by the formula
\begin{equation}
f(z) = z \kern .5pt e^{u(z) + i v(z)}.
\end{equation}
Thus a solution to (1), provided it also produces the harmonic conjugate
$v(z)$,
solves the conformal mapping problem.
Note that $u(z) + \log |z|$ is the Green's function of $\Omega$
with respect to the point~$z=0$, so $f$ is essentially the
exponential of the Green's function.
Domains of higher connectivity can also be reduced to Laplace problems.

The availability of so many tools for numerical conformal mapping
may suggest that the problem is easy, but in fact there are
challenges, of which two stand out.  One is that most domains of
practical interest contain corners, where the mapping function
will usually be singular, and it is essential to treat these
specially if one wants more than a digit or two of accuracy.
In the case of a polygon, the corners define the whole problem, and
these are dealt with by the Schwarz-Christoffel formula and its
numerical realization e.g.\ with compound Gauss-Jacobi quadrature.
The other great challenge is that of exponential distortions, which
are referred to as
the phenomenon of ``crowding''.  Whenever a domain is elongated in
certain directions, even as mildly as in the example of Figure~1,
its conformal map onto a nonelongated
canonical region will involve
exponentially large distortions.  As summarized in
Theorems 2--5 of [\rep], the distortion scales as $\exp(\pi L)$,
where $L$ is the aspect ratio of the elongation.  As a consequence,
it is usually not a good idea to attempt to
map an elongated region onto, say, a disk or a half-plane.
Other targets such as rectangles or
infinite strips may come into play, for example to treat
a domain like that of Figure~4.

I have saved the most philosophical question for last.
What is the use of numerical conformal mapping?
The following views are personal, and not all experts would
agree with them.

A great use of numerical conformal maps is to give us insight into
principles of complex analysis, harmonic functions, and their applications.
For example, the blue curves in
Figure 4 can be interpreted as flow lines of 
electricity, heat, ideal fluid, or probability
from one end to the other.  If the image rectangle has length-to-width
ratio $L$, for example, then a channel cut into this shape from a piece of
metal will have electrical resistance $\rho L$, where $\rho$ is the
resistance of a unit square.  A good image, which will almost always
have to be numerically computed, can fix these ideas beautifully in
the mind.
Throughout my career I have drawn pleasure and insight from pictures like
these.  I cannot imagine teaching complex variables without
showing some online demonstrations of conformal maps.

Specifically, two of the features that numerical conformal maps
illustrate compellingly are precisely the two
computational challenges mentioned above: behavior near singularities
(note how the blue curves in Figure 4
avoid salient corners while wrapping tightly
around reentrant ones), and exponential distortions (note the big white
region in Figure 1).

The use of conformal maps that is mentioned perhaps
more often is that they may be helpful for solving problems.
For example, every
complex analysis text tells the reader that a conformal map
may be used to solve the Laplace equation, since it reduces
a hard problem to an easy one.
I believe that the truth is not so simple.
In fact, computing a conformal map is essentially the same as solving
a Laplace problem, and whatever numerical method one employs to find
the map could probably be applied to the Laplace problem directly.
So in many cases, nothing is gained numerically from conformal mapping.
The lesson becomes even stronger for applications to other PDE\kern .5pt s that
may not be conformally invariant.

For example, suppose one is given a Laplace Dirichlet problem on the
domain on the left side of Figure 1.  One can solve it by mapping
to the disk and applying the Poisson integral formula---but where
does one get that map?  Probably by solving an integral equation or expanding
in a series, and these techniques would work equally well for the Laplace
problem itself.  In smooth multiply connected domains, for example,
series expansions work beautifully for solving
Laplace problems [\series], and although a circular
annulus is a natural domain for doubly-connected geometries, the canonical
domains with connectivity 
$\ge 3$ are less natural and do not often lead to an easy solution of your
PDE.\ \ In more extreme cases, conformal maps may not merely transplant
the difficulty but increase it, when corners or
elongations are present.

So for me, the glory of numerical conformal mapping is not in the
numbers it produces, but in the images and insights.
To finish with a fine image provided by Toby Driscoll, 
Figure 6 shows Emma's maze again, but now with the solution
indicated by a color spectrum blending from red at one end to blue at
the other.  The lightness of the colors is scaled by
$1-y^2$ if the target rectangle has its smaller dimension $-1<y<1$.
As a mark of exponential distortion or ``crowding,''
a particle beginning Brownian motion at the center of this royal road
(marked by a black dot)
would have exponentially small probability very close to
$(8/\pi) \exp(-\mu \pi/2) \approx
9.692555\times 10^{-13}$ of hitting the boundary first at one of the ends
rather than along the sides [\bornemann, chapter 10].

\begin{figure}[h]
\vspace*{7pt}
\noindent\kern 27pt\includegraphics [trim=0 0 0 0, clip, width=2.4in]{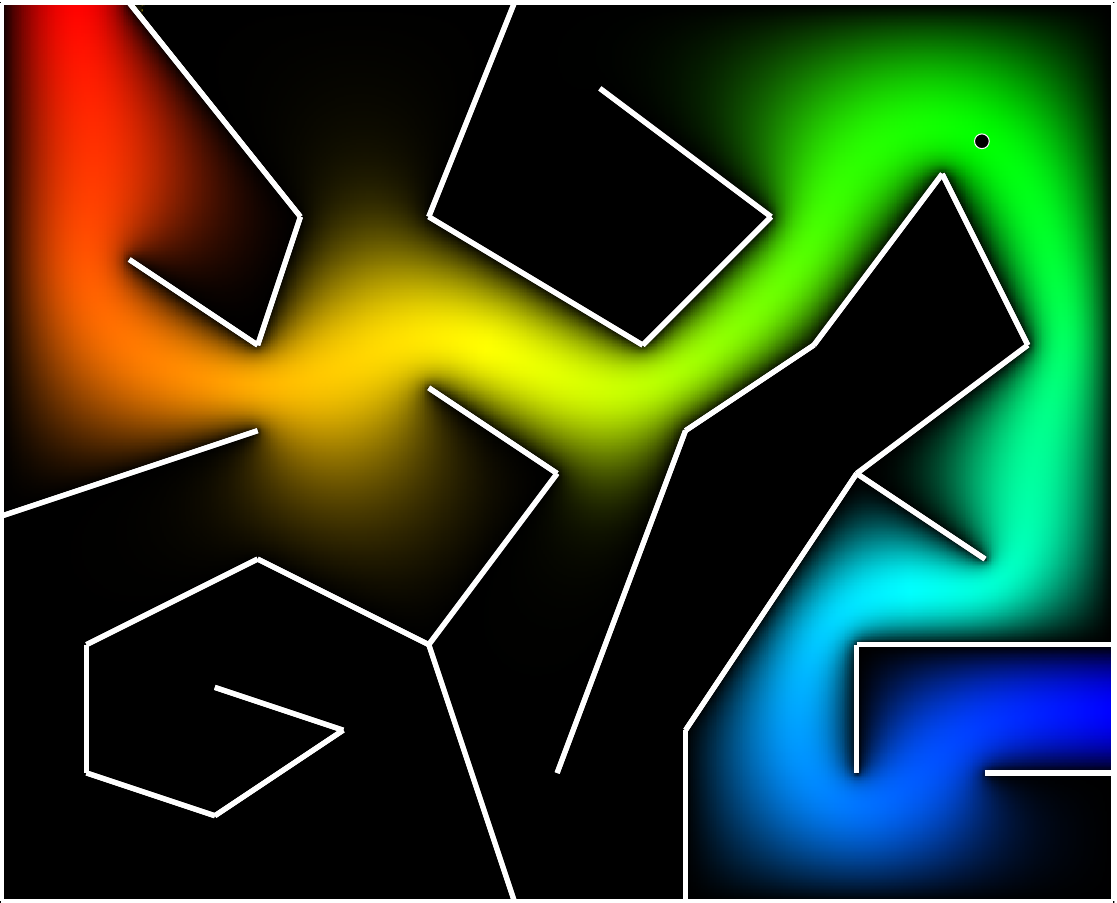}
\caption{Repetition of Figure 4 with a color spectrum showing
the conformal map---or equivalently, the
solution to a Laplace problem.}
\end{figure}

\medskip
{\bf Acknowledgments.}  
A number of colleagues have helped me considerably
in preparing this column, and I am particularly grateful
to Tom DeLillo and Toby Driscoll.

\medskip


\parindent=0pt \parskip=2pt

{\bf References}

[\amano]
K. Amano, A charge simulation method for the numerical conformal
mapping of interior, exterior and doubly-connected domains, {\em
J. Comput.\ Appl.\ Math.,} 53 (1994), 353--370.

[\banjai]
L. Banjai and L. N. Trefethen, A multipole method for
Schwarz-Christoffel mapping of polygons with thousands of sides,
{\em SIAM J. Sci.\ Comput.,} 25 (2003), 1042--1065.

[\bornemann]
F. Bornemann, D. Laurie, S. Wagon, and J. Waldvogel,
{\em The SIAM 100-digit Challenge: A Study in 
High-Accuracy Numerical Computing,} SIAM, 2004.

[\crowdy]
D. G. Crowdy, {\em Solving Problems in Multiply Connected Domains,}
SIAM (2020).

[\dekp]
T. K. DeLillo, A. R. Elcrat, E. H. Kropf, and J. A. Pfaltzgraff,
Efficient calculation of Schwarz-Christoffel transformations for
multiply connected domains using Laurent series, {\em Comput.\
Methods Funct.\ Theory,} 13 (2013), pp.~307--336.

[\dlg]
T. K. DeLillo and C. C. Green, Computation of plane potential flow
past multi-element airfoils using the Schottky-Klein prime function,
{\em Physica D,} 450 (2023), 133753.

[\driscoll] T. A. Driscoll,
Schwarz-Christoffel Toolbox in Matlab,\hfill\break
{\tt https://github.com/tobydriscoll/sc-toolbox}.

[\chebfun] T. A. Driscoll, N. Hale, and L. N. Trefethen,
{\em Chebfun Guide,} Pafnuty Publications, Oxford (2014);
see also {\tt www.chebfun.org}.

[\dt] T. A. Driscoll and L. N. Trefethen,
{\em Schwarz-Christoffel Mapping,} Cambridge (2002).

[\crdt]
T. A. Driscoll and S. A. Vavasis, Numerical conformal mapping
using cross-ratios and Delaunay triangulation,
{\em SIAM J. Sci.\ Comput.,} 19 (1998), 1783--1803.

[\gaier] D. Gaier,
{\em Konstructive Methoden der konformen Abbildung,} Springer, Berlin (1964).

[\gopal] A. Gopal and L. N. Trefethen,
Solving Laplace problems with corner singularities via
rational functions, {\em SIAM J. Numer.\ Anal.,} 57 (2019),
2074--2094; see also {\tt people.maths.ox.ac.uk/trefethen/laplace/}.

[\rep] A. Gopal and L. N. Trefethen,
Representation of conformal maps by rational
functions, {\em Numer.\ Math.,} 142 (2019), 359--382.

[\henrici]
P. Henrici, {\em Applied and Computational Complex Analysis}, v.~3,
Wiley, New York (1986).

[\zipper]
D. E. Marshall and S. Rohde, Convergence of a variant of the zipper
algorithm for conformal mapping, {\em SIAM J. Numer.\ Anal.,} 45
(2007), 2577--2609.

[\nasser]
M. M. S. Nasser, Numerical conformal mapping via a boundary integral
equation with the generalized Neumann kernel,
{\em SIAM J. Sci.\ Comput.,} 31 (2009), 1695--1715.

[\odonnell]
S. T. O'Donnell and V. Rokhlin,
A fast algorithm for the numerical evaluation of conformal mappings,
{\em SIAM J. Sci.\ Stat.\ Comput.,} 10 (1989), 475--487.

[\papa]
N. Papamichael and N. Stylianopoulos, {\em Numerical Conformal
Mapping: Domain Decomposition and the Mapping of Quadrilaterals,}
World Scientific (2010).


[\schiffer]
M. Schiffer, Some recent developments in the theory of conformal
mapping, appendix to R. Courant, {\em Dirichlet's Principle,
Conformal Mapping, and Minimal Surfaces,} Interscience (1950).

[\packing]
K. Stephenson, {\em Introduction to Circle Packing: The Theory of
Discrete Analytic Functions,} Cambridge (2005).

[\sisc80]
L. N. Trefethen, Numerical computation of the Schwarz-Christoffel
transformation, {\em SIAM J. Sci.\ Comput.,} 1 (1980), 82--102.

[\series]
L. N. Trefethen, Series solution of Laplace problems, {\em ANZIAM
J.}, 60 (2018), pp.~1--26.

[\ncm]
L. N. Trefethen, Numerical conformal mapping with rational functions,
{\em Comput.\ Methods Funct.\ Theory}, 20 (2020), pp.~369--387.

[\wegmann]
R. Wegmann, Methods for numerical conformal mapping, in {\em Handbook
of Complex Analysis: Geometric Function Theory, v.~$2$,} R.~K\"uhnau,
ed., Elsevier, 2 (2005), pp.~351--477.

\end{document}